\documentstyle[12pt]{article}
\input{amssym.def}
\input{amssym.tex}

\newcommand{\calA}{{\cal A}}

\newcommand{\qed}{{$\dashv$}}

\def\newtheorems{\newtheorem{theorem}{Theorem}[section]

                 \newtheorem{lemma}[theorem]{Lemma}
                 \newtheorem{claim}[theorem]{Claim}

                 }

\newtheorems

\title{Lusin sequences under CH and under Martin's Axiom
\thanks{The work was carried in the meeting in Lyon in 1996.}}
\author
{Uri Abraham\\
Department of Mathematics and Computer Science, \\
Ben-Gurion University, Be\'{e}r-Sheva, Israel\\
Saharon Shelah 
\thanks{Publication \# 537.}
\\
Institute of Mathematics\\
The Hebrew University, Jerusalem, Israel}
\begin{document}
\bibliographystyle{plain}
\maketitle
\begin{abstract}
Assuming the continuum hypothesis there is an inseparable sequence 
 of length $\omega_1$ that
contains no Lusin subsequence, while if Martin's Axiom
and $\neg CH$ is assumed then every inseparable sequence (of length
$\omega_1$) is a  union of
countably many Lusin subsequences.
\end{abstract}

\section{Preface}
$\omega$ is the set of natural numbers.
For $A,B \subseteq \omega$ we write $A \subseteq^* B$ iff $A\setminus B$ is
finite, and $A\perp B$ iff $A\cap B$ is finite (almost inclusion, almost
disjointness). Let $\calA=\langle A_\zeta \mid \zeta \in \omega_1 \rangle$ be
a sequence of pairwise almost disjoint subsets of $\omega_1$. So $A_\zeta
\subset \omega$ and $A_{\zeta_1}\perp A_{\zeta_1}$ for $\zeta_1 \not =
\zeta_2$. We say that $\calA$ is {\em inseparable} if there is no $B\subset
\omega$ for which
\[ (\exists^{\aleph_1} A \in \calA) (A \subseteq^* B)  \ \&\
(\exists^{\aleph_1} A \in \calA)(A \subseteq^* \omega \setminus B).\]
Thus $\calA$ is inseparable iff the collection of $B\subseteq \omega$ for
which there is an unbounded set of $\zeta$'s with $A_\zeta \subseteq^* B$ 
forms a filter on $\omega$.

   An inseparable family can be constructed in $ZFC$ alone (Lusin
\cite{1}, cited by \cite{2}). We say that $\calA$ is a {\em Lusin} sequence
if for every $i<\omega_1$ and $n\in \omega$
\[ \{ j < i \mid A_i \cap A_j \subseteq n\}\ \mbox{is finite}.\]
A seemingly stronger property is the following. We say that $\calA$ is a
{\em Lusin}* family if for every $i<\omega_1$ and $n\in \omega$
\[ \{ j<i:\; \mid A_i\cap A_j \mid < n \}\]
is finite.

It is not difficult to prove that every Lusin sequence is inseparable, and
one can construct a Lusin sequence in {\em ZFC}. Is this the only way to
build inseparable families? The answer depends on set-theoretical
assumptions as the following two results show (obtained by the first and
second author respectively).
\begin{theorem}
(1) $CH$ implies that there is an inseparable family which contains no
Lusin subsequence. (2) "Martin's Axiom$ + \neg CH$" implies that every
inseparable sequence is a countable union of Lusin* sequences.
\end{theorem}

\section{Proofs}
\subsection{$CH$ gives an inseparable non-Lusin sequence}
Assume the continuum hypothesis ($CH$) throughout this subsection.
We shall define a sequence $\calA=\langle A_\alpha \mid \alpha \in \omega_1
\rangle$ of almost disjoint subsets of $\omega$ which is inseparable
by virtue of the following property. 

\begin{equation}
\label{Pr1}
\begin{minipage}[t]{120mm}
For every infinite $X\subseteq \omega$ one of the following three
possibilities holds:
\begin{enumerate}
\item $X$ is finitely covered by $\calA$ (which means that for some finite
set $u\subset \omega_1$ $X\subseteq^* \bigcup \{ A_\alpha \mid \alpha \in u
\}$.
\item $\omega_1 \setminus X$ is finitely covered by $\calA$.
\item For some $\alpha_0<\omega_1$ for all $\alpha_0\leq \alpha < \omega_1$
$X$ splits $A_\alpha$ (which means that both $X\cap A_\alpha$ and 
$A_\alpha \setminus X$ are infinite).
\end{enumerate}
\end{minipage}
\end{equation}
It is quite obvious that if $\calA$ satisfies this property (\ref{Pr1})
then it is inseparable, and  so 
we describe the construction assuming $CH$ of a sequence that satisfies
this property, but does not contain any Lusin subsequence.

Let $\langle X_\xi \mid \xi \in \omega_1 \rangle$ be an enumeration of all
infinite subsets of $\omega$, and let $\langle e_i \mid i \in
\omega_1\rangle$ be an enumeration of all countable subsets of $\omega_1$
of order-type a limit ordinal.
The sequence $\calA=\langle A_\alpha \mid \alpha \in \omega_1\rangle$ is
defined by induction on $\alpha$. Each $A_\alpha$ is required to satisfy the
following three properties.
\begin{equation}
\label{Three}
\end{equation}
\begin{enumerate}
\item $A_\alpha\subseteq \omega$ is infinite and $A_\beta \perp A_\alpha$
for all $\beta<\alpha$.
\item For every $\xi< \alpha$ one of the following holds:
\begin{enumerate}
\item  $X_\xi$ is finitely covered by
$\langle A_\beta\mid \beta<\alpha \rangle$, or
 \item $\omega \setminus X_\xi$ is
finitely covered by $\langle A_\beta \mid \beta<\alpha \rangle$, or else
\item  $X_\xi$ splits $A_\alpha$ (both $X_\xi\cap A_\alpha$ and $A_\alpha
\setminus X_\xi$ are infinite).
\end{enumerate}
\item For every $i<\alpha$ such that $e_i \subset \alpha$ there are two
possibilities:
\begin{enumerate}
\item For some $m\in \omega$ $A_\alpha \cap A_\xi\subseteq m+1$ for an
infinite number of indices $\xi \in e_i$. (This is the ``good''
possibility.)
\item For some $m\in A_\alpha$ there is $\xi_0<\sup(e_i)$ such that for
every $\xi_0<\xi \in e_i$
\[ \min(A_\xi\setminus m+1) < \min(A_\alpha \setminus m+1).\]
(This possibility shows that $\alpha$ is not in the structure of which
$e_i$ is an elementary submodel.)
\end{enumerate}
\end{enumerate}

If we succeed then $\calA$ is inseparable (by 2) pairwise almost disjoint
sequence (by 1) which contains no Lusin subsequence (by 3) as we show.
Suppose that $I\subseteq \omega_1$ is uncountable, and $L=\langle A_i \mid
i\in I \rangle$ is a subsequence of $\calA$. We want to find some $\alpha
\in I$ and $m\in \omega$ for which $\{ \xi \in I \cap \alpha \mid A_\alpha
\cap A_\xi \subseteq m \}$ is infinite.

Consider the structure on $I$ with predicates $\omega_1$, and the binary
relation $m\in A_i$ (for $m\in \omega$, $i\in I$). Let $e\subseteq I$
be the universe of a countable elementary substructure. Then $e=e_i$ for
some $i\in \omega_1$. Let $\alpha\in I$ be any ordinal such that $\alpha >
i$ and $\alpha > \sup(e)$. We want to prove that possibility 3(a) holds for
$\alpha$. Suppose on the contrary that 3(b) holds. Then there are $m\in
A_\alpha$ and $\xi_0<\sup (e)$ such that condition 3(b) holds. As $e_i$ is
an elementary substructure, we actually have: for every $\xi_0 <\xi \in I$
 $\min(A_\xi\setminus
m+1)<\min(A_\alpha \setminus m+1)$. But this is clearly impossible for
$\xi=\alpha$ itself!

So now we turn towards the construction with its three properties listed in
\ref{Three}. To
construct $A_\alpha$ it is convenient to define a poset $P$ and a countable
collection of dense subsets of $P$. Then we define a filter $G\subseteq P$
such that $G$ intersects each of the dense sets in the countable
collection. With this we define $A_\alpha = \bigcup \{ a \mid \exists E\;
(a,E)\in G \}$, and $A_\alpha$ should have all three properties because of
the choice of the dense sets. In this fashion one does not have to
overspecify the construction.

 A condition $p=(a,E)\in P$ consists of:
\begin{enumerate}
\item A finite $a\subseteq \omega$ (which will grow to become $A_\alpha$).
\item A finite set $E\subseteq \alpha$ ($p$ promises that any extension
$a'$ of $a$ will add no new elements to $ a'\cap A_\beta$ for $\beta\in E$).
\end{enumerate}
The partial order on $P$ is defined accordingly by
\[ (a_1,E_1)\geq (a_0,E_0)\ \mbox{iff } 
\begin{minipage}[t]{100mm}
$a_1$ is an end-extension of $a_1$
such that for every $k\in a_1\setminus a_0$ and every 
$\xi \in E_0,\ k \not \in
A_\xi$ 
\end{minipage}
\]
We say that an end-extension $a_1$ of $a_0$
``respects $E$'' (where $E\subseteq \omega_1$ is finite) if $a_1\setminus
a_0 \cap A_\beta=\emptyset$ for every $\beta\in E$. So $(a_1,E_1)$ extends
$(a_0,E_0)$ iff $a_1$ is an end-extension of $a_0$ that respects $E_0$ and
$E_1$ extends $E_0$ in any way.

Now we shall define a countable collection of dense sets.
First, to ensure that $A_\alpha$ is infinite, for
every $k\in \omega$ and $p\in P$  observe that there
is an extension  $(a',E')$ with $k<\sup(a')$. Then to ensure that $A_\alpha
\perp A_\beta$ for all $\beta<\alpha$ observe that $(a,E\cup\{\beta\})$
extends $p$.
This takes care of 1.

For every $X\subseteq \omega$ such that neither $X$ nor $\omega\setminus X$
are covered by $\langle A_\beta \mid \beta < \alpha \rangle$, and for every
$k\in \omega$, define $D_{X,k} \subset P$ by: $(a,E)\in D_{X,k}$ iff both
$a\cap X$ and $a\setminus X$ contain $\geq k$ members.
\begin{claim}
$D_{X,k}$ is dense (open) in $P$.
\end{claim}
{\bf Proof.}
Suppose that $(a_0,E_0)\in P$.
Since neither $X$ nor its complement are $\subseteq^*$-included in
$A=\bigcup\{ A_\beta \mid \beta \in E_0 \}$, both $X\setminus A$ and
$(\omega\setminus X)\setminus A$ are infinite. We can find an end-extension
$a_1$ of $a_0$ such that $(a_1\setminus a_0)\cap A=\emptyset$ and both
$a_1\cap X$ and $a_1\setminus X$ contain $\geq k$ members.
Thus $(a_0,E_0)<(a_1,E_0)\in D_{X,k}$. \qed

So add to the countable list of dense sets all sets $D_{X_\xi,k}$ for $k\in
\omega$ and $\xi<\alpha$ such that neither $X_\xi$ nor $\omega \setminus
X_\xi$ are finitely covered by $\langle A_\beta \mid \beta < \alpha
\rangle$. This ensures 2.

The main issue of the proof is to take care of 3. What dense sets will do
the job? Fix $e=e_i$ for some $i<\alpha$ such that $e\subseteq \alpha$. We
say that $p=(a,E)\in P$ is of type (a) if for some $m\in a$ the following
holds.
\begin{equation}
\begin{minipage}[c]{110mm}
For every end-extension $a'$ of $a$ that respects $E$ ($a'\setminus a$ is
disjoint to $\bigcup \{ A_\beta \mid \beta
\in E \}$) and for every $\xi_0 \in e$ there is some $\xi\in e$, $\xi_0 \leq
\xi$, such that $A_\xi \cap a' \subseteq m+1$.
\end{minipage}
\end{equation}

If $p$ is of type (a) then
\begin{enumerate}
\item the least $m\in a$ that satisfies  is denoted $m_p$,
\item any extension of $p$ is also of type (a) (with the same $m$).
\end{enumerate}

We say that $p=(a,E)\in P$ is of type (b) if there are two successive
members of $a$, $m$ and $n$ (so $m,n\in a$ and $(m,n)\cap a = \emptyset$).
such that for some $\xi_0\in e$ for every $\xi_0 \leq \xi \in e$ $A_\xi
\cap (m,n) \not = \emptyset$.
\begin{claim}
Any condition in $P$ has an extension of type (a) or an extension of type
(b).
\end{claim}
{\bf Proof.}
Given $p=(a,E)$ let $m=\max(a)$. Is $p$ of type (a) by virtue of $m$? If
yes, we are done, and if not then there are
\begin{enumerate}
\item  $a'$ an end-extension of $a$,
respecting $ E$, and
\item $\xi_0 \in e$,
\end{enumerate}
such that for every $\xi \in e$ with $\xi_0 \leq \xi$ $A_\xi \cap a'
\setminus m+1\not = \emptyset$.
Let $n> \max a'$ be such that $n\not \in \bigcup \{ A_\beta \mid \beta \in
E \}$ and consider the condition $p'=(a\cup\{ n \},E)$ extending $p$. Then
for every $\xi_0 \leq \xi \in e$ $A_\xi \cap (m,n)\not = \emptyset$. That
is, $p'$ is of type (b).\qed

For every $\xi_0\in e$ define $D_{\xi_0,e}$ by 
$p=(a,E)\in D_{\xi_0,e}$ iff either $p$
is of type (b) or $p$ is of type (a) and there exists 
some $\xi \in e\cap E$
above $\xi_0$ with  $A_\xi\cap a \subseteq m+1$ (where
$m=m_p$). $D_{\xi_0,e}$ is dense by the following argument.
 Suppose $p_0\in P$ is
given. If $p_0$ is extendible into a condition of type (b) then we are
done. Assume that $p_1=(a_1,E_1)
 \geq p_0$ is of type (a). By definition of type (a),
there is some $\xi \in e$, with $\xi_0 \leq \xi$ such that $A_\xi \cap a_1
\subseteq m+1$. Hence $(a_1, E_1\cup \{ \xi \})\in D_{\xi_0,e}$ as
required. 

Add to the countable list of dense sets all sets
$D_{\xi_0,e}$ where $e=e_i$ for some $i<\alpha$ such that $e_i \subseteq
\alpha$ and $\xi_0\in e_i$. This ensures property 3 (of \ref{Three}).
Suppose that $A_\alpha$ is defined from the filter $G$ that intersects all
the above dense sets. Given $i$ such that $e_i=e \subseteq \alpha$, we ask:
is there $(a,E)\in G$ of type (b) for $e$? If yes, then possibility 3(b)
holds for $A_\alpha$. If no, then there is an $m$ and an unbounded set of
$\xi\in e$ such that $A_\xi \cap A_\alpha\subseteq m+1$. To see this,
observe that $m_p=m$ is fixed for all $p\in G$ of type (a). Now consider
any $\xi_0\in e$, and pick $p=(a,E)\in D_{\xi_0,e}\cap G$. Then $p$ is of
type (a) and there is
$\xi \in E\cap e_i$ above $\xi_0$ with $A_\xi \cap a \subseteq m+1$. Thus
$A_\xi \cap A_\alpha \subseteq m+1$ follows.

\subsection{Martin's Axiom: Inseparable $\Rightarrow$ contains a Lusin*
subsequence}
Assume Martin's Axiom $+ 2^{\aleph_0}>\aleph_0$.
Let $\calA=\langle A_\zeta \mid \zeta \in \omega_1 \rangle$ be an
inseparable sequence of length $\omega_1$ (any length below the continuum
works). Define the following poset:
\[ Q = \{ (u,n)\mid u\subseteq \omega_1\ \mbox{is finite and } n < \omega
\}\]
ordered by
\begin{equation}
 (u_1,n_1) \leq (u_2,n_2)\ \mbox{iff }  
\begin{minipage}[t]{75mm}
 $u_1 \subseteq u_2\ \& \ n_1\leq n_2\ \&$ \\
 $\forall i\in u_1\; \forall j\in u_2\setminus u_1\ j<i \Rightarrow \mid
A_i \cap A_j \mid > n_1.$
\end{minipage}
\end{equation}
This is easily shown to be transitive. We intend to prove that $Q$ is a
c.c.c poset and that for every $\alpha <\omega_1$ and $k<\omega$ the 
set $D_{\alpha,k}$ of
$(u,n)$ in $Q$ for which $\sup(u)>\alpha$ and $n>k$ is dense. So if
$G\subset Q$ is a filter provided by Martin's Axiom which intersects each
of these dense sets, then $U=\bigcup \{ u \mid \exists n (u,n)\in G \}$ is
uncountable and $\langle A_\alpha \mid \alpha \in U \rangle$ is a Lusin*
sequence. (Because if $i\in U$ and $k<\omega$ then  the collection of $A_j$
for $j \in U\cap i$ such that $\mid A_i\cap A_j\mid \leq k$ is finite by the
following argument.  For some $(u,n)\in G$ $i\in u$ and $n\geq k$. This
implies that $\mid A_i \cap A_j \mid > k$ for every $j<i$ such that $j\in
U\setminus u$.)

The full result, concerning the decomposition of $\calA$ into countably many
Lusin* subsequences, follows from the fact that (under Martin's Axiom)
if $Q$ is a c.c.c poset, then $Q$ is a countable union of filters (each
intersects the required dense sets). (Consider the poset $Q^{<\omega}$
consisting of finite sequences from $Q$, ordered coordinate-wise.)

It is easy to see that if $(u,n)\in Q$ and $v$ is any end-extension of $u$
then $(u,n)\leq (v,n)$. Also, if $n \leq m$
then $(u,n)\leq (u,m)$. This shows that the required sets $D_{\alpha,k}$ 
are indeed dense
in $Q$, and so the main point of the proof is to show that $Q$ satisfies the
c.c. condition.

\begin{lemma} $Q$ satisfies the c.c. condition.
\end{lemma}
{\bf Proof}.
Let $\langle (u_\zeta,n_\zeta)\mid \zeta \in \omega_1 \rangle$ be an
$\omega_1$-sequence of conditions in $Q$. We may assume that for some fixed
$n$ and $k$ $n=n_\zeta$  and $k=\mid u_\zeta \mid$ for all
$\zeta \in \omega_1$, and that the sets $u_\zeta's$ form a $\Delta$ system.
That is, for some finite  $c_0\subset \omega_1$ $c_0 = u_{\zeta_1}\cap
u_{\zeta_2}$ for all $\zeta_1 \not = \zeta_2$ and $\max(u_{\zeta_1}) <
\min(u_{\zeta_2} \setminus c_0)$ for $\zeta_1 < \zeta_2$. 

We want to find
$\zeta_1 < \zeta_2$ such that $(u_{\zeta_1}\cup u_{\zeta_2},n)$ extends
both $(u_{\zeta_1},n)$ and $(u_{\zeta_2},n)$.
It is evident that $(u_{\zeta_1}\cup u_{\zeta_2},n)$ extends
$(u_{\zeta_1},n)$ (the lower part) but the problem is the possibility that
for some $i\in u_{\zeta_2}$ and $j \in u_{\zeta_1}\setminus c_0$ $\mid A_i
\cap A_j \mid \leq n$.

We shall find two uncountable sets $K,L\subseteq \omega_1$ such that for
every $\zeta_1\in K$ and $\zeta_2\in L$ $(u_{\zeta_1},n)$ and
$(u_{\zeta_2},n)$ are compatible. We start with $K_0=L_0=\omega_1$, and
define $K_{i+1}\subseteq K_i$ and $L_{i+1}\subseteq L_i$ by induction for
$i<\mid u_\zeta \setminus c_0\mid$ (any $\zeta$ as these sets have all the
same size). The definition of $K_i$ and $L_i$ depends on a finite
parameter set, and it is convenient to have a countable model in which the
definition is carried on. So let $M\prec \langle H_{\omega_1},\calA, Q,
\{(u_\zeta,n_\zeta): \zeta \in \omega_1\}\rangle$ be a countable elementary
submodel. The following lemma is used.
\begin{lemma}
Let $U,V\in M$ be two uncountable subsets of $\omega_1$ and $n<\omega$.
There are uncountable subsets $U_1\subseteq U$ and $V_1\subseteq V$
(definable in $M$) such that for every $\zeta \in U_1$ and $\xi \in V_1$
$\mid A_\zeta \cap A_\xi\mid > n$ (and hence $(\{\zeta\},n)$ and
$(\{\xi\},n)$ are compatible in $Q$).
\end{lemma}
It should be obvious how successive applications of the lemma yield the
c.c.c., and so we turn to the proof of the lemma. Let $\delta =\omega_1\cap
M$ be the set of countable ordinals in our countable structure $M$.

\noindent
{\bf Case 1}: for some $\zeta \in U \setminus \delta$ and $\xi \in V
\setminus \delta$ $\mid A_\zeta \cap A_\xi\mid > n$.
In this case pick $X \subset A_\zeta \cap A_\xi$ with $\mid X \mid >n$,
and let $U_1=\{ i \in U \mid X \subset A_i \}$, $V_1=\{j\in V \mid X
\subset A_j \}$. Both $U_1$ and $V_1$ are uncountable (for if $U_1$ is
countable then it would be included in $M$, but $A_\zeta$ shows that this is
not the case).

\noindent
{\bf Case 2}: not Case 1. So for every $\zeta \in U \setminus \delta$ and
$\xi \in V
\setminus \delta$ $\mid A_\zeta \cap A_\xi\mid\leq n$. Let $0\leq m_0 \leq
n$ be the maximal size of some intersection $F=A_\zeta \cap A_\xi$ for
indices $\zeta$ and $\xi$ as above. Then
$U_1=\{ i \in U \mid F \subset A_i\}$ and $V_1 = \{ j\in V \mid X \subset
A_j\}$ are uncountable and for $i\in U_1$ and $j\in V_1$ $A_i\cap A_j=F$
(by maximality of $\mid F \mid$). So the sets $A_\zeta \setminus F$ and
$A_\xi$ for  $\zeta \in U \setminus \delta$ and $\xi \in V
\setminus \delta$ are pairwise disjoint, which gives a contradictory
separation of $\calA$.

\end{document}